\documentclass[12pt,thmsa]{article}
\usepackage{sw20bams}

\input tcilatex
\begin{document}

\author{Michael Schweitzer and Steven Finch}
\title{Zero Divisors in Associative Algebras over Infinite Fields}
\date{March 30, 1999}
\maketitle

\begin{abstract}
Let $F$ be an infinite field. We prove that the right zero divisors of a
three-dimensional associative $F$-algebra $A$ must form the union of at most
finitely many linear subspaces of $A$. The proof is elementary and written
with students as the intended audience.
\end{abstract}

\section{Introductory Example}

Consider three-dimensional real space $R^3$ endowed with the following
vector multiplication: 
\[
\left( 
\begin{array}{c}
\alpha \\ 
\beta \\ 
\gamma
\end{array}
\right) \cdot \left( 
\begin{array}{c}
\delta \\ 
\epsilon \\ 
\varphi
\end{array}
\right) =\left( 
\begin{array}{c}
\alpha \delta +\frac 12(\gamma \varphi -\beta \epsilon )-\frac{\sqrt{3}}%
2(\beta \varphi +\gamma \epsilon ) \\ 
(\alpha \epsilon +\beta \delta )-\frac{\sqrt{2}}8(5\gamma \varphi -\beta
\epsilon )-\frac{\sqrt{6}}8(\beta \varphi +\gamma \epsilon ) \\ 
(\alpha \varphi +\gamma \delta )+\frac{\sqrt{6}}8(\gamma \varphi +3\beta
\epsilon )-\frac{\sqrt{2}}8(\beta \varphi +\gamma \epsilon )
\end{array}
\right) 
\]
This definition first appeared in \cite{Lucas} and has application to
pattern recognition and information processing.

Observe that multiplication is associative, commutative and has identity
element $\left( 1,0,0\right) $. To further understand this algebra, it is
natural to study the set $Z$ of all zero divisors, that is, all vectors
which are mapped to zero under multiplication by a nonzero element. A
determinant argument gives that Z is the set of all $\left( \alpha ,\beta
,\gamma \right) $ satisfying the homogeneous cubic equation 
\[
6\sqrt{6}\gamma ^3+(6\sqrt{2}\beta -12\alpha )\gamma ^2+(6\sqrt{6}\beta ^2+24%
\sqrt{3}\alpha \beta )\gamma -10\sqrt{2}\beta ^3+12\alpha \beta ^2+16\alpha
^3=0 
\]
It may surprise some that this equation factors: 
\[
\left( \sqrt{6}\gamma +2\alpha -\sqrt{2}\beta \right) \left( \left( 2\alpha -%
\sqrt{2}\beta -\sqrt{6}\gamma \right) ^2+\left( 2\alpha +2\sqrt{2}\beta
\right) ^2\right) =0 
\]
and hence $Z$ is the union of a plane and a line passing through the origin.
In fact, this is not surprising to algebraists \cite{Atiyah}, \cite{Zariski}%
. For any finite-dimensional, associative, commutative algebra $A$ with
identity, the set $Z$ is the union of at most finitely many linear subspaces
of $A$. We say $A$ is \textbf{tame} when $Z$ is so simply described. Does
tameness hold even if the commutativity hypothesis is discarded? No, the
four-dimensional algebra of real $2\times 2$ matrices is not tame. What can
be said about the three-dimensional case? We answer this question in the
following sections. As far as is known, this material has not previously
appeared in the literature.

\section{Main Result}

Let $F$ be an infinite field, for example, the real numbers $R$ or the
complex numbers $C$. Let $A$ be a finite-dimensional $F$-algebra, that is, $A
$ is a finite-dimensional vector space over $F$ together with a bilinear map 
$A\times A\longrightarrow A$, known as vector multiplication $xy$, which
satisfies $f(xy)=(fx)y=x(fy)$ for all $f\in F$ and $x,y\in A$. Further, let $%
A$ be associative, that is, $(xy)z=x(yz)$ for all $x,y,z\in A$. We do not
assume $A$ is commutative or that $A$\ has an identity element. Everything
in the following holds true if one systematically interchanges ''left'' and
''right''.

A subring $I$ of $A$ is a \textbf{left ideal} if $a\in A$, $x\in I$ implies
that $ax\in I$. A left ideal $I\neq A$ is \textbf{maximal} if, whenever $J$
is a left ideal such that $I\subset J\subset A$, then either $J=I$ or $J=A$.

An element $z\in A$ is a \textbf{right zero divisor} if there exists a
nonzero $w\in A$ with $wz=0$. The set of all right zero divisors is denoted
by $Z$. Note that we consider $0\in Z$, which is uncustomary.

The algebra $A$ is \textbf{right tame} if $Z$ is the union of at most
finitely many linear subspaces of $A$. We believe this phraseology to be
new. Define $A$ to be \textbf{right proper tame} if $A$ is right tame and,
additionally, $Z\neq A$. Also define $AA$ to be the set of all products $xy$%
, where $x\in A$, $y\in A$. Our main result depends on whether $AA$ fills
out all of $A$ or not; the condition $AA=A$ may be thought of as a poor
man's replacement for the existence of an identity element.\\

\noindent \textbf{Theorem. }(i)\textit{\ If }$AA\neq A$\textit{, then }$Z=A$%
\textit{\ and hence }$A$\textit{\ is right tame.}

(ii) \textit{If }$AA=A$\textit{, then }$Z$\textit{\ is the union of all
maximal left ideals of }$A$\textit{. Further, }$A$\textit{\ is right proper
tame if and only if }$A$\textit{\ possesses at most finitely many maximal
left ideals.}\\

\noindent \textbf{Corollary.} \textit{If }$dim(A)=3$\textit{, then A is
right tame.}\\

The proofs of these assertions and preliminary lemmas are given in the next
section. For the benefit of students and for completeness' sake, we provide
many underlying details. The required background is essentially covered in 
\cite{Herstein}.

\section{Detailed Proofs}

Under the conditions on $A$ described in section 2, if $S$ is a subset of $A$%
, define 
\begin{eqnarray*}
L(S) &=&\text{ the linear subspace generated by }S \\
&=&\text{ the intersection of all linear subspaces which contain }S
\end{eqnarray*}
and 
\begin{eqnarray*}
I(S) &=&\text{ the left ideal generated by }S \\
\ &=&\text{ the intersection of all left ideals which contain }S\text{.}
\end{eqnarray*}
In the following, if $k$ is a positive integer and $s\in S$, then $ks$ is
repeated addition: $ks=\tsum\limits_{j=1}^ks$, $0s=s$ and $(-k)s=-(ks)$.

\begin{lemma}
$L(S)=\left\{ \sum\limits_{i=1}^nf_is_i:f_i\in F,s_i\in S,\text{any integer }%
n\geq 1\right\} $ and

$I(S)=\left\{ \sum\limits_{i=1}^na_is_i+\sum\limits_{j=1}^mk_jt_j:a_i\in
A,s_i,t_j\in S,\text{any integers }k_j,\text{ }n\geq 1,m\geq 1\right\} $.
\end{lemma}

\noindent \textbf{Proof. }Focus only on the second formula. First, show that
the right-hand side is a left ideal containing $S$. It clearly contains $S$
since $1$ is an integer and $1s=s$ for any $s\in S$. It is a left ideal
since, given $b\in A$, 
\[
b\cdot \left( \sum\limits_ia_is_i+\sum\limits_jk_jt_j\right)
=\sum\limits_ib(a_is_i)+\sum\limits_jb(k_jt_j)=\sum\limits_i(ba_i)s_i+\sum%
\limits_j(k_jb)t_j 
\]
by associativity and since both $ba_i\in A$ and $k_jb\in A$. Conversely,
show that the right-hand side is contained in every left ideal containing $S$%
.\\

\begin{lemma}
If $I$ is a left ideal, then $L(I)$ is also a left ideal and $A\cdot L(I)$
is contained in $I$.
\end{lemma}

\noindent \textbf{Proof. }If $f_i\in F,x_i\in I$ and $a\in A$, then\textbf{\ 
}$a\cdot \left( \sum\limits_if_ix_i\right) =\sum\limits_if_i(ax_i)\in L(I)$
since $ax_i\in I$, so $L(I)$ is a left ideal. Also, $\sum\limits_if_i(ax_i)=%
\sum\limits_i(f_ia)x_i\in I$, so $A\cdot L(I)\subset I$.\\

Let $x\in A$ and define $\left\langle x\right\rangle =I(\left\{ x\right\} )$%
, the left ideal generated by $\left\{ x\right\} $. Define too a map $%
R(x):A\longrightarrow A$ by $R(x)(y)=yx$, the right-multiplication map by $x$%
.\\

\begin{lemma}
The map $R(x)$ is a linear transformation, $Im(R(x))\subset \left\langle
x\right\rangle $ and $A\cdot \left\langle x\right\rangle \subset Im(R(x))$.
In particular, if $\left\langle x\right\rangle =A$, then $Im(R(x))=AA.$
\end{lemma}

\noindent \textbf{Proof. }The first two assertions are clear. By Lemma 1, $%
\left\langle x\right\rangle $ is the set of all $ax+kx$ where $a\in A$ and $%
k $ is an integer, so for any $b\in A$, 
\[
b\cdot \left( ax+kx\right) =(ba+kb)\cdot x=R(x)(ba+kb) 
\]
by associativity, hence $A\cdot \left\langle x\right\rangle \subset Im(R(x))$%
. If $\left\langle x\right\rangle =A$, then $AA=A\cdot \left\langle
x\right\rangle $ is contained in $Im(R(x))$, but $Im(R(x))$ is trivially
contained in $AA$.\\

\begin{lemma}
Let $x\in A$. Then $x$ is a right zero divisor if and only if $Im(R(x))\neq
A $.
\end{lemma}

\noindent \textbf{Proof. }$x$ is a right zero divisor iff there exists $%
y\neq 0$ with $R(x)(y)=0$. This is possible iff $Ker(R(x))\neq \left\{
0\right\} $. This, in turn, is possible iff $Im(R(x))\neq A$.\\

\noindent \textbf{Proof of Theorem, part (i). }We prove, assuming $AA\neq A$%
, that every element of $A$ is both a right and a left zero divisor. Let $%
x\in A$, $x\neq 0$. Then $Im(R(x))\subset AA\neq A$, which implies
immediately that $x$ is a right zero divisor by Lemma 4. Likewise, $x$ is a
left zero divisor.\\

Having dealt with the easy case for which $AA\neq A$, we devote attention to
the more interesting case $AA=A$. Given two sets $S$ and $T$, define $S+T$
to be the set of all sums $s+t$, where $s\in S$, $t\in T$.\\

\begin{lemma}
Assume $AA=A$. If $I$ is a maximal left ideal, then $I=L(I)$. That is, a
maximal left ideal is necessarily a linear subspace. If $I$ and $J$ are
maximal left ideals and $I\neq J$, then $I+J=$ $A$.
\end{lemma}

\noindent \textbf{Proof. }By Lemma 2, $L(I)$ is a left ideal and $I\subset
L(I)\subset A$. By maximality, it follows that $I=L(I)$ or $L(I)=A$. The
latter is impossible since otherwise $A=AA=A\cdot L(I)\subset I$ by Lemma 2,
but $I\neq A$ by definition. Thus $I=L(I)$. Also, $I+J$ is a left ideal, $%
I\subset I+J$, and $I+J\neq I$ (since otherwise $J=\left\{ 0\right\}
+J\subset I$, contradicting maximality of $J$). Therefore $I+J$ contains $I$
properly and, by maximality of $I$, this implies that $I+J=A$.\\

\begin{lemma}
Let $V$ be a vector space over $F$. Then $V$ is not the union of finitely
many proper linear subspaces.
\end{lemma}

\noindent \textbf{Proof. }Suppose not. Among all possible such
decompositions of V, choose one 
\[
V=\bigcup\limits_{i=1}^nL_i\text{ } 
\]
with $n$ minimal. Then $n\geq 2$ and no $L_i$ is contained in the union of
the other $L_j$'s. Select vectors $v_1$ and $v_2$ such that $v_1\in L_1$ but 
$v_1\notin L_i$ for $i\neq 1$, and $v_2\in L_2$ but $v_2\notin L_i$ for $%
i\neq 2$. Consider the affine line $K=\left\{ v_1+fv_2:f\in F\right\} $ and
suppose that it intersects $L_i$ in two distinct points $p=v_1+fv_2$ and $%
q=v_1+gv_2$, $f\neq g$. Then $p-q=(f-g)\cdot v_2$ and $gp-fq=(g-f)\cdot v_1$
are also in $L_i$, but this is impossible since $v_1$ is in $L_1$ only and $%
v_2$ is in $L_2$ only. So $K$ intersects each $L_i$ in at most one point.
Since V is the union of $L_1,L_2,...,L_n$, it follows that $K$ is a finite
set. This is absurd since $F$ is infinite. \\

\noindent \textbf{Proof of Theorem, part (ii). }We first prove, assuming $%
AA=A$, that the set $Z$ of all right zero divisors is the union of all
maximal left ideals of $A$. Let $x\in Z$, then $Im(R(x))\neq A$ by Lemma 4.
If $\left\langle x\right\rangle =A$, then $Im(R(x))=AA=A$ by Lemma 3, which
is a contradiction. Hence $\left\langle x\right\rangle \neq A$ and $%
\left\langle x\right\rangle $ must be contained in some maximal left ideal.
Conversely, let $x$ be an element of a maximal left ideal. Thus $%
\left\langle x\right\rangle \neq A$. It follows from Lemma 3 that $%
Im(R(x))\subset \left\langle x\right\rangle $, hence $Im(R(x))\neq A$.
Therefore by Lemma 4, $x\in Z$.

If there are only finitely many maximal left ideals, then by Lemmas 5 and 6,
each of them is a linear subspace and thus $A$ is right proper tame.
Conversely, if $A$ is right proper tame, then there exist linear subspaces $%
L_1,L_2,...,L_n$ such that 
\[
Z=\bigcup\limits_{i=1}^nL_i\text{ = the (possibly infinite) union of maximal
left ideals.} 
\]
Let $M$ be a maximal left ideal and let $M_i=M\bigcap L_i$. Then 
\[
M=\bigcup\limits_{i=1}^nM_i\text{ } 
\]
By Lemma 5, $M$ is itself a vector space. A finite decomposition of $M$ as
linear subspaces $M_i$ is impossible, by Lemma 6, unless $M=M_i$ for some $i$%
. Thus $M\subset L_i$. We wish to show that $L_i$ contains no other maximal
left ideal. Suppose $N\subset L_i$ for some maximal left ideal $N\neq M$.
Then $A=M+N\subset L_i$ by Lemma 5 and hence $Z=A$, which cannot be true.
Therefore every maximal left ideal is contained in some $L_i$ and no two are
contained in the same $L_i$, so they must be only finite in number.\\

To prove the Corollary, we need some properties of polynomials in $n$
variables with coefficients in $F$. An ideal $P$ in the ring $F\left[
X_1,X_2,X_3,...,X_n\right] $ (in fact, any ring) is \textbf{prime} if for
all $p,q\in F\left[ X_1,X_2,X_3,...,X_n\right] $, $p\cdot q\in P$ implies $%
p\in P$ or $q\in P$. \\

\begin{lemma}
If a polynomial $p\in $ $F\left[ X_1,X_2,X_3,...,X_n\right] $ vanishes over
all of $F^n$, then $p=0$.
\end{lemma}

\noindent \textbf{Proof. }If $n=1$ and $p\neq 0$, then $p$ can have at most
finitely many zeros, but $F$ is infinite. Use induction on $n$ to complete
the argument.\\

\begin{lemma}
If $p,q\in $ $F\left[ X_1,X_2,X_3,...,X_n\right] $, $p$ is linear and every
zero of $p$ is also a zero of $q$, then $p$ divides $q$.
\end{lemma}

\noindent \textbf{Proof. }Make an invertible affine change of coordinates so
that $p$ becomes $X_1$. Then the zero set of $p$ is the set of all $%
(0,f_2,..,f_n)$, $f_i\in F$. Write 
\[
q=q_0(X_2,..,X_n)+q_1(X_2,..,X_n)\cdot X1+... 
\]
then 
\[
0=q(0,f_2,..,f_n)=q_0(f_2,...,f_n) 
\]
for all $f_2,...,f_n\in F$. By Lemma 7, $q_0=0$ and hence $X_1$ divides $q$.%
\\

\begin{lemma}
Assume $p,q,r\in $ $F\left[ X_1,X_2,X_3,...,X_n\right] $ and $p$ and $q$ are
both linear. If $p$ and $q$ both divide $r$ and $q$ is not a scalar multiple
of $p$, then $p\cdot q$ divides $r$.
\end{lemma}

\noindent \textbf{Proof. }Again, change coordinates so that $p$ becomes $X_1$%
. By assumption, $r=p\cdot u=q\cdot v$ for some polynomials $u$ and $v$,
hence $X1\cdot u=q\cdot v$. The ideal generated by $X1$ is prime since the
quotient is $F\left[ X_2,X_3,...,X_n\right] $, an integral domain \cite
{Hungerford}. We deduce $q$ is not in this ideal since $q$ is linear but is
not a scalar multiple of $X_1$. Thus $v=X_1\cdot w$ for some polynomial $w$.
Therefore $r=q\cdot X_1\cdot w$ and so $X_1\cdot q$ divides $r$.\\

The above propositions can be strengthened, of course, but these are all we
need.

Consider now the \textbf{determinant form} $D(\xi _1,...,\xi _n)$ of the
algebra $A$, that is, the determinant of right multiplication $R(x)$
relative to a fixed basis $b_1,...,b_n$. We have $x=\sum\limits_i\xi
_ib_i\in Z$ iff $D(\xi _1,...,\xi _n)=0$ by Lemma 4. From Lemma 7, clearly $%
Z=A$ iff $D=0$. In the following, assume $AA=A$ as always and $dim(A)=n$. A
linear subspace of $A$ is said to have \textbf{codimension} $k$ if it has
dimension $n-k$. The guiding principle is that it is easy to bound the
number of maximal left ideals of either low dimension or low codimension.\\

\begin{lemma}
The algebra $A$ possesses at most $n$ maximal left ideals of codimension 1.
\end{lemma}

\noindent \textbf{Proof. }Let $M$ be a maximal left ideal of codimension 1. $%
M$ is the zero set of a linear form $f$. Every element of $M$ is a right
zero divisor by the Theorem, hence a zero of $D$. By Lemma 8, $f$ divides $D$%
. If $N$ is another maximal left ideal of codimension 1, it is the zero set
of some linear form $g$. Since $M\neq N$, $g$ is not a scalar multiple of $f$
and $g$ also divides $D$. By Lemma 9, $f\cdot g$ divides $D$. Since the
degree of D is $\leq n$, there can be at most $n$ maximal left ideals of
codimension 1.\\

\begin{lemma}
The algebra $A$ possesses at most one maximal left ideal of dimension $<n/2$.
\end{lemma}

\noindent \textbf{Proof. }If $M$ and $N$ are maximal left ideals of
dimension $<n/2$ and $M\neq N$, then on one hand $M+N=A$ by Lemma 5. On the
other hand, $\dim (M+N)\leq \dim (M)+\dim (N)<n$, which is a contradiction.
Hence $M=N$.\\

\begin{lemma}
If $A$ possesses a maximal left ideal of dimension 1, then $A$ has at most $%
n+1$ maximal left ideals. In particular, $A$ is right tame.
\end{lemma}

\noindent \textbf{Proof. }Let $M$ be that ideal. If $N$ is some other
maximal left ideal, then $M+N=A$ by Lemma 5. Since 
\[
n=\dim (M+N)\leq \dim (M)+\dim (N)=1+\dim (N) 
\]
it follows that $\dim (N)\geq n-1.$ Since $N$ is a proper subspace (being a
maximal ideal), $\dim (N)=n-1.$ Hence $N$ has codimension 1. Now use Lemma
10.\\

\noindent \textbf{Proof of Corollary. }A maximal left ideal has dimension 0,
1 or 2. If there is one with dimension 0, then it is the only one by Lemma
11. If there is one with dimension 1, then there are at most four by Lemma
12. Otherwise all are of dimension 2 and there are at most three by Lemma 10.

\section{Closing Remarks}

We have proved that a three-dimensional associative $F$-algebra is right
tame, where $F$ is infinite. For the one or two-dimensional cases, the
associativity requirement may be dropped and right tameness still follows.
The case when $F$ is finite is also trivial.

Using some tools from algebraic geometry \cite{Hartshorne}, we can prove the
following general result:\\

\noindent \textbf{Theorem. }\textit{Let }$A$\textit{\ be a
finite-dimensional algebra over an algebraically closed field }$F$\textit{.
Let }$D$\textit{\ be the determinant form for right multiplication with
respect to some basis of }$A$\textit{. Then A is right tame iff }$D$\textit{%
\ splits into a product of linear forms over }$F$\textit{.}\\

\noindent \textbf{Proof. }Assuming $A$ is right tame and $D\neq 0$, let $f$
be an irreducible factor of $D$. The ideal $\left( f\right) $ generated by $%
f $ is prime, thus the zero set $V$ of $f$ is a variety. Since $V$ is
contained in a finite union of linear subspaces, it must be contained in one
of them. Call this subspace $L$. Under the correspondence between algebraic
sets and ideals of polynomial rings, the ideal $I(L)$ is contained in the
ideal $I(V)$, but $I(V)=\left( f\right) $ because $\left( f\right) $ is
prime. $I(V)$ is generated by linear forms, hence these must be multiples of 
$f$. Therefore $f$ is linear.\\

We ask the following:\\

\noindent \textbf{Open Question. }\textit{Let }$A$\textit{\ be a
finite-dimensional algebra over }$R$\textit{\ and let }$D$\textit{\ be as
before. Is A right tame iff }$D$\textit{\ splits into a product of linear
forms over }$C$ ?\\

For an infinite-dimensional algebra, there is in general no close
relationship between zero divisors and maximal ideals. For example, let $A$
be the group algebra for the group of (rational) integers over a field $F$.
Then A has no zero divisors ($\neq 0$) but it has maximal ideals $\neq 0$.
So in the infinite-dimensional case, one can assume the existence of an
identity element and commutativity, etc., but all this doesn't help.

Finally, returning to the three-dimensional real case, there is
computational evidence that the associativity requirement may be weakened
somewhat \cite{Finch} and yet right tameness still holds. An successful
elementary treatment (as above) under such extended circumstances is not
likely.

\end{document}